\documentclass[11pt]{article}
\usepackage{amssymb}
\usepackage{a4wide}
\usepackage{pslatex}
\usepackage{theorem}
\usepackage{fancyhdr}

\DeclareMathAlphabet{\mathcal}{OMS}{cmsy}{m}{n}

\title{A note on the geometry of linear Hamiltonian systems \\ of signature 0 in $\RR^4$}
\author{James Montaldi}

\def\eqref#1{(\ref{#1})}

\def\lin{{\Gamma}}
\def\CC{\mathbb{C}}
\def\HH{\mathbb{H}}
\def\KK{\mathbb{K}}
\def\RR{\mathbb{R}}
\def\ZZ{\mathbb{Z}}
\def\nullcone{\mathcal{N}}
\def\GL{\mathbf{GL}}
\def\JJ{\mathbf{J}}
\def\OO{\mathbf{O}}
\def\SO{\mathbf{SO}}
\def\UU{\mathbf{U}}
\def\Sp{\mathbf{Sp}}
\def\dd{\mathrm{d}}
\def\Hom{\mathop{\mathrm{Hom}}\nolimits}
\def\tr{\mathop{\mathrm{tr}}\nolimits}
\def\DD{\mathcal{D}}

\def\half{{\textstyle\frac12}}

\newtheorem{theorem}{Theorem}
\newtheorem{lemma}[theorem]{Lemma}
\newtheorem{proposition}[theorem]{Proposition}

{\theorembodyfont{\normalfont}
\newtheorem{example}[theorem]{Example}

}

\newenvironment{proof}{\noindent\textit{Proof}\ }{\hfill
  $\Box$\bigskip}

\begin{document}

\maketitle

\begin{abstract}
\noindent It is shown that a linear Hamiltonian system of signature
zero on $\RR^4$ is elliptic or hyperbolic according to the number of
Lagrangian planes in the null-cone $H^{-1}(0)$, or equivalently the
number of invariant Lagrangian planes. An extension to higher
dimensions is described.
\end{abstract}

\section*{Introduction}

A linear Hamiltonian system is determined by a quadratic form
$H$---the Hamiltonian---and a symplectic form $\omega$.  In terms of
a basis, let $S$ be the symmetric matrix representing $H$ (so $H(x)
= x^TSx$) and $\Omega$ the nondegenerate skew-symmetric matrix
representing $\omega$ (so $\omega(x,y)=x^T\Omega y$). The linear
dynamical system is then given by $\dot x=Lx$ where
$L=\Omega^{-1}S$.

The Hamiltonian system $(H,\omega)$ or $(S,\Omega)$ is
\emph{degenerate} if $S$ is degenerate, otherwise it is
\emph{nondegenerate}. A nondegenerate system is said to be
\emph{elliptic} if all eigenvalues of $L$ are pure imaginary,
\emph{hyperbolic} if no eigenvalues of $L$ are pure imaginary, and
\emph{mixed} if it is neither elliptic nor hyperbolic.

Given a linear Hamiltonian system $(H,\omega)$ it is in principle
straightforward to determine its type (elliptic/hyperbolic/ mixed)
simply by calculating the eigenvalues of $L$. On the other hand, the
type depends only on the symplectic geometry of the quadratic form
$H$, and this elementary note is a first attempt to make this
dependence explicit. For a system in $\RR^4$, we relate the type to
the number of Lagrangian planes in the \emph{null-cone}
$\nullcone=H^{-1}(0)$ of $H$.  The main result is:

\begin{theorem} \label{thm:main}
Let $(\RR^4, \omega, H)$ be a nondegenerate linear Hamiltonian
system, with $H$ of signature 0, and let $\nullcone$ be the
null-cone of $H$. Then
\begin{itemize}
\item if $H$ is elliptic with simple eigenvalues,  $\nullcone$ contains no
Lagrangian planes; %
\item if $H$ is elliptic with double eigenvalues and non-zero nilpotent part,
 $\nullcone$ contains precisely 1 Lagrangian plane; %
\item if $H$ is hyperbolic with a full quadruplet of eigenvalues,
$\nullcone$ contains precisely 2 Lagrangian planes; %
\item if $H$ is hyperbolic with a pair of coincident real eigenvalues and
non-zero nilpotent part,  $\nullcone$ contains precisely 3
Lagrangian planes; %
\item if $H$ is hyperbolic with simple real eigenvalues, $\nullcone$
contains precisely 4 Lagrangian planes.
\end{itemize}
\end{theorem}

In the course of the proof, we also show that if the Hamiltonian
system has double eigenvalues and is semisimple, then the null-cone
contains infinitely many Lagrangian planes (and ``$\infty+2$'' if
the double eigenvalues are real).

It should perhaps be pointed out that signature 0 is the interesting
case. If $H$ has signature $\pm2$ (so is positive or negative
definite) then the system is elliptic, and $H^{-1}(0)=\{0\}$ which
contains no planes at all, Lagrangian or otherwise. And if $H$ is of
signature $\pm1$, then it is necessarily of mixed type and again
$\nullcone$ contains no planes (Lagrangian or otherwise). In fact,
as we see below, a non-degenerate quadratic form on $\RR^{2n}$ whose
null-cone contains a subspace of dimension $n$ is necessarily of
vanishing signature.

In higher dimensions, there is no such one-to-one correspondence
between the type of equilibrium and the number of Lagrangian planes
in the null-cone of a signature zero quadratic Hamiltonian.

\section{Linear geometry of the null-cone of a quadratic form}
\label{sec:null-cone}

Although we do not need the full generality, we discuss arbitrary
non-degenerate quadratic forms in this section as the arguments are
sufficiently simple and short.  Let $H:\RR^n\to\RR$ be a quadratic
form of index $\ell$, and denote its null-cone $H^{-1}(0)$ by
$\nullcone$. We assume that $2\ell\leq n$, as otherwise we can
replace $H$ by $-H$, which both have the same null-cone. After a
suitable change of basis we can write $H:\RR^\ell\oplus\RR^k\to\RR$
(where $k=n-\ell$) so that $H(x,y)=\|y\|^2 - \|x\|^2$.  Let
$\nullcone=H^{-1}(0)$, and let $\lin\subset \nullcone$ be a linear
space. Let $\Pi:\RR^\ell\oplus\RR^k\to\RR^\ell$ be the Cartesian
projection on to the first factor. Note that since $\lin \subset
\nullcone$ we know that $\Pi|_\lin$ is injective, for otherwise
$\lin\cap\RR^k\neq\{0\}$ which would contradict the obvious fact
that $\nullcone\cap\RR^k=\{0\}$. It follows that any linear subspace
of the null-cone has dimension at most $\ell$. Furthermore, since
$\lin\cap\RR^k=\{0\}$ it follows that $\lin$ can be represented as
the graph of a linear map $\RR^\ell\to\RR^k$, with matrix $M$. That
is, $\lin=\lin_M=\{(x,y)\in\RR^\ell\oplus\RR^k \mid y=Mx\}$.

Suppose now $\lin=\lin_M\subset\nullcone$. Then $H(Mx,x) = \|Mx\|^2
- \|x\|^2 \equiv 0$ (for all $x\in\RR^\ell$). It follows that the
$\ell$ columns of $M$ are orthonormal vectors in $\RR^k$. In
particular, if $k=\ell$, $M$ is an orthogonal matrix. Consequently,

\begin{proposition} \label{prop:frames}
Let $H$ be a quadratic form in $n$ variables of index $\ell\leq
n/2$.  Then the null-cone of $H$ contains no linear spaces of
dimension greater than $\ell$, and the set of linear spaces of
dimension $\ell$ contained in the null-cones is in 1-1
correspondence with the set of orthonormal $\ell$-frames in
$\RR^{n-\ell}$. In the case that $H$ is of signature zero, so
$n=2\ell$, this set can be identified with the group of orthogonal
$\ell\times\ell$ matrices.
\end{proposition}

The set of orthonormal $\ell$-frames in $\RR^n$ is called a Stiefel
manifold, denoted $V_{\ell,n}$.  If $\ell<n$ then $V_{\ell,n}$ is
connected, while if $\ell=n$ it has two connected components,
distinguished by the determinant of $M$.

In particular, in the case $n=4$ and $\ell=2$, the set of planes in
$\nullcone$ can be identified with the union of two circles. Indeed,
if $(x,y)\in\RR^2\oplus\RR^2$ is such that $\|y\|^2=\|x\|^2$ then
there are precisely two orthogonal matrices $M,M'$ with $y=Mx=M'x$,
and $\det(M)=1=-\det(M')$. Conversely, given any $M,M'\in\OO(2)$
with $\det(M)=1=-\det(M')$ then $\lin_M\cap\lin_{M'}$ consists of a
single line through the origin. The two circles therefore lie in
different connected components of the Stiefel manifold $V_{2,4}$.

Finally, recall that the Grassmannian $G_{2,4}$ of planes in $\RR^4$
is of dimension 4, while the Lagrange-Grassmannian $\Lambda_2$ is of
dimension 3.  We are interested in the intersection between
$\Lambda_2$ and the 1-dimensional family of planes in $\nullcone$
inside $G_{2,4}$ which could generically be expected to be finite.

\section{Invariant Lagrangian planes}
Any Lagrangian plane in $\nullcone$ is invariant under the
Hamiltonian dynamics, so that this ``symplectic geometry'' of the
null-cone is intimately related to the dynamics:

\begin{lemma} \label{lem:invariant}
  Let $(M,\omega, H)$ be any (smooth) Hamiltonian system, and let
  $L$ be a Lagrangian submanifold of $M$. Then $L$ is invariant if
  and only if it is contained in a level set of the Hamiltonian.
\end{lemma}

\begin{proof}
$L$ is invariant if and only if it is tangent to the characteristic
direction of the Hamiltonian system at each of its points.  Since
$L$ is Lagrangian, this is equivalent to the tangent space to $L$ at
each point annihilating $\dd H$, which in turn is equivalent to $L$
being contained in a level set of $H$. (This argument continues to
hold even at critical points of $H$.)
\end{proof}

\begin{lemma}\label{lem:decompose}
Suppose $(V,\omega,H)$ is a linear Hamiltonian system, which
decomposes into the direct sum of two subsystems $(V_1,\omega_1,
H_1) \oplus (V_2,\omega_2, H_2)$, and let $\nullcone_j$ be the
null-cone of $H_j$. Suppose that there are no eigenvalues common to
the two subsystems $H_1$ and $H_2$. Then every Lagrangian subspace
$L \subset \nullcone$ can be decomposed as a direct sum $L=L_1\oplus
L_2$ such that each $L_j$ is a Lagrangian subspace of\/ $V_j$
contained in $\nullcone_j$.
\end{lemma}

\begin{proof}
This involves the lemma above that Lagrangian subspaces in the
null-cone are invariant.  In general, any subspace invariant under a
linear transformation is a direct sum of (generalized) eigenspaces.
Since any (generalized) eigenspace of $(V,\omega,H)$ is contained in
either $V_1$ or $V_2$, it follows that $L=L_1\oplus L_2$ with
$L_j\subset V_j$ and invariant. That each $L_j$ is Lagrangian in
$V_j$ follows easily by contradiction, and since $L_j$ is Lagrangian
in $V_j$ and invariant it is contained in $\nullcone_j$, by Lemma
\ref{lem:invariant}.
\end{proof}

\section{Lagrangian planes in the null-cone of a Hamiltonian} \label{sec:proof}

In this section we prove Theorem \ref{thm:main}, proceeding case by
case. The argument in each case is either based on Lemma
\ref{lem:decompose} (if the system is a product of 2-dimensional
systems) or on a calculation involving a normal form for the
Hamiltonian (see \cite{MMCM}). Given a quadratic form $H$ on $\RR^4$
of index 2, we know from Section \ref{sec:null-cone} that the
null-cone $\nullcone$ contains 2 circles of planes. Rather than use
a basis adapted to the quadratic form, we use a symplectic basis in
order to make the Lagrangian nature of a given plane transparent.
Recall that if a plane is given by $y=Mx$ and the symplectic form is
$\omega=\dd y_1\wedge \dd x_1 + \dd y_2\wedge \dd x_2$, then the
plane is Lagrangian if and only if $M$ is symmetric (in which case
$M$ is the hessian matrix of the generating function determining the
plane).

The proof in the hyperbolic and non-semisimple elliptic cases
consist in taking the Hamiltonian $H(q,p)$ in some normal form,
substituting $y=Mx$ for an arbitrary symmetric matrix $M$, where $y$
and $x$ are suitable (symplectic) choices of $p$'s and $q$'s, and
finally determining for which such $M$ the restriction of $H$ to the
graph of $y=Mx$ vanishes. In order to consider all possible
Lagrangian planes it is necessary to consider 4 cases:
(i) $y=(p_1,p_2)$ and $x=(q_1,q_2)$; %
(ii) $y=(q_1,q_2)$ and $x=(-p_1,-p_2)$; %
(iii) $y=(q_1,p_2)$ and $x=(-p_1,q_2)$; %
and (iv) $y=(p_1,q_2)$ and $x=(q_1,-p_2)$ %
(the signs are chosen so that the symplectic form $\dd p_1\wedge \dd
q_1 + \dd p_2\wedge \dd q_2 = \dd y_1\wedge \dd x_1 + \dd y_2\wedge
\dd x_2$). In the calculations we will always take
$M=\pmatrix{\alpha&\beta\cr \beta&\gamma}$, often without further
reference.

\subsection{The elliptic cases: imaginary eigenvalues}

\paragraph{The semisimple system}
 If $H$ has distinct imaginary eigenvalues $\pm
i\lambda_1$ and $\pm i\lambda_2$, ($\lambda_1,\lambda_2$ positive
and distinct) the result follows from Lemma \ref{lem:decompose}, for
each of the two ``modes'' (eigenspaces for $\pm i\lambda_1$ and for
$\pm i \lambda_2$) is symplectic and has no invariant Lagrangian
subspaces.

If on the other hand, $\lambda_1=\lambda_2$ and the system is
semisimple, then every non-zero point is contained in a periodic
orbit, and that periodic orbit spans a plane in $\RR^4$ which
consists entirely of periodic orbits, so is invariant. So if the
initial point lies in $\nullcone$, the orbit is in $\nullcone$ and
so therefore is the plane it spans. That this 2-dimensional plane is
Lagrangian follows from Lemma \ref{lem:invariant}.

Consequently, if $\lambda_1=\lambda_2$ every point in $\nullcone$ is
contained in a Lagrangian plane in $\nullcone$, and so there are
infinitely many such planes (forming one of the two families of
planes in $\nullcone$~; the other consists of symplectic planes).

\paragraph{The non-semisimple elliptic system}
Generically, an elliptic Hamiltonian with double eigenvalues has a
nontrivial nilpotent part.  For a normal form we take
$$
H_\pm = \pm\half(q_1^2+q_2^2) +\lambda(p_2q_1 - p_1q_2)
$$
for which the associated linear system has eigenvalues $\pm
i\lambda$, with multiplicity 2 (we assume $\lambda\neq0$).  The two
normal forms $H_\pm$ are symplectically inequivalent.

We show that the only Lagrangian plane in the null-cone is
$\{q_1=q_2=0\}$. There are 4 cases to consider as was pointed out
above.

\noindent\textbf{(i)} $y=(p_1,p_2)$ and $x=(q_1,q_2)$. With
$M=\pmatrix{\alpha&\beta\cr \beta&\gamma}$ we substitute $y=Mx$ into
the normal form for the Hamiltonian, to obtain
$$
2H_\pm = (\pm1+2\lambda\beta)q_1^2 + 2\lambda(\gamma-\alpha)q_1q_2
+ (\pm1-2\lambda\beta)q_2^2
$$
It is clear that this cannot vanish identically, so there are no
Lagrangian planes in this portion of $\nullcone$.

\noindent\textbf{(ii)} $y=(q_1,q_2)$ and $x=(-p_1,-p_2)$. With
$M=\pmatrix{\alpha&\beta\cr \beta&\gamma}$ we again substitute
$y=Mx$ into the normal form for the Hamiltonian, to obtain the more
complicated expression
$$
2H_\pm = (\pm(\alpha^2+\beta^2) + 2\lambda\beta)p_1^2 +
2(\pm\beta(\alpha+\gamma) + \lambda(\gamma-\alpha))p_1p_2
+(\pm(\beta^2+\delta^2)-2\lambda\beta)p_2^2
$$

Suppose this vanishes identically; then if $\beta=0$ it follows that
$\alpha=\gamma=0$, while if $\beta\neq0$ the coefficients of $p_1^2$
and $p_2^2$ cannot both vanish.  The only Lagrangian plane in
$\nullcone$ therefore corresponds to $M=0$ which is therefore
$\{q_1=q_2=0\}$.

\noindent\textbf{(iii)} $y=(q_1,p_2)$ and $x=(-p_1,q_2)$, and (iv)
$y=(p_1,q_2)$ and $x=(q_1,-p_2)$.  Similar computations to those
above show that there are no Lagrangian planes of this form.

The single Lagrangian plane found in this part is the limit of the
(symplectic) normal modes in the semisimple case above, as the
eigenvalues approach equality, and indeed on this plane the motion
is periodic with period $2\pi/\lambda$.

\subsection{The hyperbolic case: a quadruplet of eigenvalues}
Here we take
$$H = \kappa(p_1q_1+p_2q_2) + \lambda(p_1q_2 - p_2q_1).$$
The eigenvalues of the linear system are $\pm\kappa\pm i\lambda$. We
consider the 4 cases as before:

\noindent\textbf{(i)} $y=(p_1,p_2)$ and $x=(q_1,q_2)$. %
Substituting $y=Mx$, with  as usual $M=\pmatrix{\alpha&\beta\cr
\beta&\gamma}$, into the Hamiltonian gives
$$
 H = (\kappa\alpha-\lambda\beta)q_1^2 +
 (2\kappa\beta+\lambda(\alpha-\gamma))q_1q_2 + (\kappa\gamma+\lambda\beta)q_2^2.
$$
This only vanishes identically if $\alpha=\beta=\gamma=0$,
corresponding to the plane $\{p_1=p_2=0\}$.

\noindent\textbf{(ii)} $y=(q_1,q_2)$ and $x=(-p_1,-p_2)$. %
Substituting this $y=Mx$ (same $M$ as before) into the Hamiltonian
gives
$$
 H = (\kappa\alpha+\lambda\beta)p_1^2 +
 (2\kappa\beta+\lambda(\gamma-\alpha))p_1p_2 + (\kappa\gamma-\lambda\beta)p_2^2.
$$
Again, this only vanishes identically if $\alpha=\beta=\gamma=0$,
corresponding to the plane $\{q_1=q_2=0\}$.

\noindent\textbf{(iii)} $y=(q_1,p_2)$ and $x=(-p_1,q_2)$ and (iv)
$y=(p_1,q_2)$ and $x=(q_1,-p_2)$. Repeating similar calculations
shows in these cases that there are no other solutions.

Thus the null-cone contains precisely 2 Lagrangian planes.  These
two planes are in fact the stable and unstable manifolds of the
vector field: every initial point in $\{q_1=q_2=0\}$ tends to the
origin as $t\to\infty$ (the stable manifold), while every initial
point in $\{p_1=p_2=0\}$ tends to the origin as $t\to-\infty$ (the
unstable manifold).

Further calculations show that the two Lagrangian planes belong to
the same family of planes in $\nullcone$. Indeed, write
$H=p_1(\kappa q_1+\lambda q_2) + p_2(-\lambda q_1+\kappa q_2) \equiv
p_1 Q_1 + p_2Q_2$. Thus,
 $$4H=(p_1+Q_1)^2-(p_1-Q_1)^2 + (p_2+Q_2)^2-(p_2-Q_2)^2.$$
In these coordinates, $\pmatrix{p_1+Q_1\cr p_2+Q_2} =
M\pmatrix{p_1-Q_1\cr p_2-Q_2}$ lies in $\nullcone$ iff $M\in\OO(2)$.
The two Lagrangian planes found above correspond to $M=-I$ and $M=I$
respectively, and both are in $\SO(2)$.

\subsection{The hyperbolic cases: real eigenvalues}

\paragraph{Semisimple cases}
Here the normal form is
$$H= \lambda_1p_1q_1 + \lambda_2p_2q_2,$$
and the associated linear system has eigenvalues $\pm\lambda_1,
\pm\lambda_2$, and we can assume $\lambda_1>0, \lambda_2>0$. In the
2-dimensional hyperbolic system with $H=\lambda pq$ there are two
invariant (Lagrangian) lines.  It then follows from Lemma
\ref{lem:decompose} that if $\lambda_1\neq\lambda_2 $ there are
precisely 4 invariant Lagrangian planes in the 4-dimensional system,
as required. We now show this again by direct calculation, as the
calculation is required for the case $\lambda_1=\lambda_2$.

Taking the four cases for Lagrangian planes in turn:

\noindent\textbf{(i)} $y=(p_1,p_2)$ and $x=(q_1,q_2)$. %
With $M=\pmatrix{\alpha&\beta\cr \beta&\gamma}$ as usual, we
substitute $y=Mx$ into the Hamiltonian to obtain
$$H =\lambda_1\alpha q_1^2 + (\lambda_1+\lambda_2)\beta q_1q_2 +
\lambda_2\gamma q_2^2.
$$
Since $\lambda_1,\lambda_2>0$, this vanishes if and only if $M=0$,
corresponding to the plane $\{p_1=p_2=0\}$.

\noindent\textbf{(ii)} $y=(q_1,q_2)$ and $x=(-p_1,-p_2)$. %
Again, since $\lambda_1,\lambda_2>0$, the Hamiltonian restricted to
$q=Mp$ vanishes if and only if $M=0$, corresponding to the plane
$\{q_1=q_2=0\}$.

\noindent\textbf{(iii)} $y=(q_1,p_2)$ and $x=(-p_1,q_2)$. %
Substituting $y=Mx$ gives
$$H =\lambda_1\alpha p_1^2 + (\lambda_2-\lambda_1)\beta p_1q_2 -\lambda_2\gamma
q_2^2.
$$
Provided $\lambda_1\neq\lambda_2$, this vanishes again only when $M=0$,
corresponding to the plane $\{q_1=p_2=0\}$.

\noindent\textbf{(iv)} $y=(p_1,q_2)$ and $x=(q_1,-p_2)$. %
Again, provided $\lambda_1\neq\lambda_2$, the restriction of the Hamiltonian
vanishes again only when $M=0$, corresponding to the plane $\{p_1=q_2=0\}$.

If $\lambda_1=\lambda_2$, so the system has a double real eigenvalue
and is semisimple, then in cases  (i) and (ii) the same conclusion
applies, while in cases (iii) and (iv) the planes $\Gamma_M$ are in
the null-cone whenever $\alpha=\gamma=0$ (for all $\beta$).  That
is, of the two circles of planes in $\nullcone$ (Section
\ref{sec:null-cone}), one consists entirely of Lagrangian planes
while the other contains exactly two Lagrangian planes.

The planes $\{p_1=p_2=0\}$ and $\{q_1=q_2=0\}$ (cases (i) and (ii))
are respectively the stable and unstable manifolds of the vector
field.  On the other hand, the remaining invariant planes are all
saddles. These last ones are not the only invariant planes with
saddles (for example $\{q_1=p_1=0\}$ is an invariant symplectic
plane with a saddle point), but they are distinguished by the fact
that the flow on a Lagrangian plane does not need to be
area-preserving.

\paragraph{Non-semisimple case}
A normal form is
$$H = \lambda(p_1q_1+p_2q_2) + p_1q_2.$$
The eigenvalues of the corresponding system are $\pm\lambda$.

\noindent\textbf{(i)} $y=(p_1,p_2)$ and $x=(q_1,q_2)$. %
Substituting $y=Mx$ as usual, one obtains
$$H = \lambda\alpha q_1^2 + (\alpha+2\lambda\beta)q_1q_2 +
(\lambda\gamma+\beta)q_2^2.
$$
This vanishes identically if and only if $M=0$, corresponding to the
plane $\{p_1=p_2=0\}$ (the unstable manifold).

\noindent\textbf{(ii)} $y=(q_1,q_2)$ and $x=(-p_1,-p_2)$. %
In this case we obtain $H=\lambda\alpha p_2^2 +
(\gamma+2\lambda\beta)p_1p_2 + (\lambda\alpha+\beta)p_1^2$, which
vanishes only when $M=0$, corresponding to the plane $\{q_1=q_2=0\}$
(the stable manifold).

\noindent\textbf{(iii)} $y=(q_1,p_2)$ and $x=(-p_1,q_2)$. %
In this case there are no solutions.

\noindent\textbf{(iv)} $y=(p_1,q_2)$ and $x=(q_1,-p_2)$. %
Here the only solution is $\{p_1=q_2=0\}$. (On this plane the
dynamics is a simple area-preserving saddle, with eigenvalues
$\pm\lambda$.)

\section{Higher dimensional systems} \label{sec:Rn}

In higher dimensions the number of Lagrangian planes in a null-cone
does not distinguish between different types of equilibrium.  For
example,  in $\RR^8$ with its standard symplectic structure, the two
Hamiltonians
$$
\begin{array}{lcl}
H_1 &=&
(p_1^2+q_1^2)-2(p_2^2+q_2^2)+3(p_3^2+q_3^2)-4(p_4^2+q_4^2)\\[6pt]
H_2 &=& (p_1^2+q_1^2)-2(p_2^2+q_2^2)+ p_3q_3 + 2 p_4q_4
\end{array}
$$
are both of signature 0. The first is elliptic while the second has
two pairs of imaginary eigenvalues and two pairs of real ones, yet
neither has any invariant Lagrangian subspaces.

The fact that neither has any Lagrangian subspaces in its null-cone
follows from the following argument.

Let $(V,\omega, H)$ be a linear Hamiltonian system with
$\mathrm{signature}(H)=0$, and suppose that the system is of
codimension at most 1.  Then $(V,\omega, H)$ splits into a direct
sum of symplectic ``eigenspaces"
$$(V,\omega, H) = \bigoplus_{j=1}^r (V_j,\omega_j, H_j),$$
such that the eigenvalue-quadruplets
$\{\lambda_j,-\lambda_j,\bar\lambda_j,-\bar\lambda_j\}$ of the
different components are distinct. The codimension hypothesis
implies that each component is of dimension $2$ or $4$.

\newpage

Define a function $\delta(j)$ by
\begin{itemize}
\item If $j$ is such that $\dim V_j=2$ and $\lambda_j$ is
imaginary then $\delta(j)=0$;

\item If $j$ is such that $\dim V_j=4$ and $\lambda_j$ is imaginary
(so of algebraic multiplicity 2 and geometric multiplicity 1), then
$\delta(j)=1$;

\item If $j$ is such that $\dim V_j=4$ and $\lambda_j$ has non-zero
real and imaginary parts, then $\delta(j)=2$;

\item If $j$ is such that $\dim V_j=4$ and $\lambda_j$ is real (so of
algebraic multiplicity 2 and geometric multiplicity 1), then
$\delta(j)=3$.

\item If $j$ is such that $\dim V_j=2$ and $\lambda_j$ is real then
$\delta(j)=2$;
\end{itemize}

Then it follows immediately from Theorem \ref{thm:main} and Lemma
\ref{lem:decompose} that the number of Lagrangian planes in the
null-cone $\nullcone$ is $\prod_{j=1}^r\delta(j)$.

\paragraph{Acknowledgements} The ideas in this note were elaborated during a visit
to Cincinnati, and I would like to thank Ken Meyer for several useful discussions.

\bigskip

\hbox to 0.35\textwidth{\hrulefill}

\smallskip

 \parindent=0pt

{\sl School of Mathematics, \\
University of Manchester, \\
Manchester M60 1QD, \\
UK.}

\texttt{j.montaldi@manchester.ac.uk}

\hbox to 0.35\textwidth{\hrulefill}

\end{document}